

\input amstex
\documentstyle{amsppt}
\magnification=1200
\catcode`\@=11
\redefine\logo@{}
\catcode`\@=13

\define \bn{\Bbb N}
\define \bz{\Bbb Z}
\define \bq{\Bbb Q}
\define \br{\Bbb R}
\define \bc{\Bbb C}
\define \bh{\Bbb H}

\define \M{{\Cal M}}

\define\La{{\Cal L}}
\define\geg{{\goth g}}
\define\0o{{\overline 0}}
\define\1o{{\overline 1}}


\define\mult{\text{mult}}
\define\re{\text{re}}
\define\im{\text{im}}

\define\NEF{\text{NEF}}
\define\ir{\text{ir}}
\define\EF{\text{EF}}
\TagsOnRight

\document

\topmatter

\title
K3 surfaces, Lorentzian Kac--Moody algebras and Mirror Symmetry
\endtitle

\author
Valeri A. Gritsenko \footnote{Supported by
SFB 170 ``Geometrie und Analysis''. \hfill\hfill} and
Viacheslav V. Nikulin \footnote{Partially supported by
Grant of Russian Fund of Fundamental Research;
Grant of ISF MI6000; Grant of ISF and
Russian Government MI6300; Grant of AMS; SFB 170
``Geometrie und Analysis''. \hfill\hfill}
\endauthor

\address
St. Petersburg Department Steklov Mathematical Institute
Fontanka 27,
\newline
${}\hskip 8pt $
191011 St. Petersburg,  Russia
\endaddress
\email
gritsenk\@pdmi.ras.ru;\ \ gritsenk\@cfgauss.uni-math.gwdg.de
\endemail

\address
Steklov Mathematical Institute,
ul. Vavilova 42, Moscow 117966, GSP-1, Russia.
\endaddress
\email
slava\@nikulin.mian.su
\endemail

\abstract
We consider the variant of Mirror Symmetry Conjecture
for K3 surfaces which
relates ``geometry" of curves of a general member
of a family of K3 with ``algebraic functions" on the moduli
of the mirror family.
Lorentzian Kac--Moody algebras are involved in this construction.
We give several examples when this conjecture is valid.
\endabstract

\rightheadtext
{K3, Kac--Moody and Mirror symmetry}

\leftheadtext{V. Gritsenko and  V. Nikulin}

\endtopmatter

\document

\subhead
0. Introduction
\endsubhead

In this paper we want to interpret our results \cite{GN} and
\cite{N10} from the view-point of the mirror symmetry for K3 surface.
This interpretation was the subject of the talk given by the second
author at the conference ``Toric geometry" in  Warwich University at
September 18--23, 1995.
\vskip2pt
This paper was written during our stay
at SFB 170  ``Geometrie und Analysis'' of Mathematical
Institute of Georg-August-University
in  September -- October 1995.
We are grateful to the Institute for hospitality.

\subhead
1. Mirror symmetry for K3 surfaces
\endsubhead

Let $S$ be an even hyperbolic lattice, i.e. a free $\bz$-module of
rank $n+1$
with an integral even symmetric bilinear form of the signature
$(1,n)$. This lattice $S$ may appear in two ways in connection
with algebraic K3 surfaces:

\smallpagebreak

(A) $S=S_X$ is a Picard lattice of a K3 surface $X$. These
K3 surfaces form a family
$$
\M_S=\{\ \text{K3 surface\ } X\ |\ S\subset S_X\}
$$
of dimension $20-\dim S$ (see \cite{N1}, \cite{N6} for
definition which is actually based on local \cite{G.N. Tyurina, \u S}, and
global Torelli
Theorem \cite{P-\u S--\u S} for K3 surfaces and epimorphicity of Torelli map
\cite{Ku} for K3 surfaces).
A general member $X$ of this family has the Picard lattice
$S_X=S$.

\smallpagebreak

(B) $S =([c]^\perp_T)/[c]$, for the lattice of transcendental cycles
$T=T_X$ (transcendental lattice) of
a K3 surface $X$ where $c \in T$ is a primitive element of $T$ with
$c^2=0$. (We consider K3 surfaces over $\bc$; then
$T_X=(S_X)^\perp_{H^2(X,\bz)}$.)
These K3 surfaces $X$ form a family
$$
\M_{T^\perp}=\{ \ \text{K3 surface\ } X \ |\ T_X \subset T\}
$$
of dimension $\dim S$. A general member $X$ of
this family has $T_X=T$.

\smallpagebreak

These two families $\M_S$ and $\M_{T^\perp}$
are called dual (or mirror symmetric, or mirror). This is how Mirror
Symmetry for K3 surfaces
(inspired by explanation of the strange Arnol'd duality \cite{A})
had first appeared in \cite{P}, \cite{DN} and \cite{N2}, \cite{D1}.
In particular, in \cite{N2}
there was developed some lattice theory for the
exact calculation of these dual families. The new understanding
of mirror symmetry for K3 (see for example \cite{D2})
which is due to the modern
Physics and the Mirror Symmetry for Calabi--Yau 3-folds
(see \cite{COGP} and \cite{Mor1}, \cite{Mor2}) is related with
the fact that one can calculate the lattice $S$ for the situation
(B) using Yukawa coupling at the point defined by $c$
at infinity of $\M_{T^\perp}$.

For the model (A), the lattice $S$ is related with the geometry of
curves and is the intersection form
of all curves on a general K3 surface $X\in \M_S$.
For the model (B), the lattice
$S$ is related with the geometry of moduli $\M_{T^\perp}$ at an
appropriate point at infinity of the mirror (dual)
family $\M_{T^\perp}\subset \overline{\M}_{T^\perp}$. Thus,
for any question related with the geometry of curves
on a general member $X$ of the family
$\M_S$, one can ask about its analog for the dual family $\M_{T^\perp}$
from the point of view of the geometry of the moduli $\M_{T^\perp}$.

An effect we want to discuss here is the following:

\smallpagebreak

{\it It turns out that in some cases ``geometry" of irreducible and
effective classes of
divisors of general
$X \in \M_S$ is related with interesting ``algebraic
functions" on the dual family $\M_{T^\perp}$. This relation
involves Lorentzian Kac--Moody Lie algebras
(and conjecturally some physics).}

\smallpagebreak

Here an element of $S=S_X$ is called {\it irreducible} (respectively
{\it effective}) if it contains an irreducible (respectively effective) curve.
The moduli $\M_{T^\perp}$ is a quotient of a symmetric domain of
type IV by some arithmetic group $G$, and ``algebraic function" means here
an automorphic form with respect to $G$ on this domain.

\subhead
2. Geometry of irreducible and effective classes of divisors on a K3 surface
\endsubhead

In this section, we consider a hyperbolic lattice $S$ from
the point of view of
the model (A). Thus, now $S=S_X$ is the Picard lattice of a K3 surface
$X$.
Then elements of $S$ reflect some geometry of curves.
An element $h \in S$ is called {\it irreducible} if it contains
an irreducible curve on $X$. An element $h \in S$ is
called {\it effective} if it is a finite sum of irreducible elements.
For K3 surfaces, effective and irreducible classes may be
described (up to automorphisms of $S$) purely arithmetically
using only the intersection form of the lattice $S$.

Now we give this description.
It is sufficient to describe the set  $\Delta^{\ir}\subset S$
of all irreducible elements.
It is well-known (and very easy to see) that $h^2 \ge -2$ if
$h \in \Delta^{\ir}$. An irreducible
element $\delta\in \Delta^{\ir}$ with $\delta ^2=-2$
contains a non-singular irreducible rational curve (exceptional curve)
on $X$. In particular,
$\Delta^{\ir} = \Delta^{\ir}_{-2}\cup  \Delta^{\ir}_{\ge 0}$ where
$$
\Delta^{\ir}_{-2} =\{ \delta \in \Delta^{\ir}\ |
\delta^2=-2\},\ \  \Delta^{\ir}_{\ge 0}=
\{ h \in \Delta^{\ir} \ |\ h^2\ge 0\}.
$$
An element $h \in S$ is called {\it nef} if
$h\cdot C\ge 0$ for any irreducible curve $C$ on $X$. We denote by
$\NEF(S)$ the set of all nef elements of $S$. It is known that
the sets $\NEF(S)$ and $\Delta^{\ir}_{\ge 0}$ almost coincide.
Obviously, $x^2\ge 0$ if $x \in \NEF(S)$ and
$\bn \Delta^{\ir}_{\ge 0}\subset \NEF (X)$. If
$c \in \NEF(S)$ and $c^2=0$, then $c \in  \Delta^{\ir}_{\ge 0}$
if and only if $c$ is primitive \cite{P-\u S--\u S}.
If $h \in \NEF(S)$ and $h^2>0$, then
$h \in  \Delta^{\ir}_{\ge 0}$ if and only if
there does not exist primitive $c \in  \NEF (S)$ with
$c^2=0$ such that  $c\cdot h=1$ (in particular,
$2h \in \Delta^{\ir}_{\ge 0}$).
See \cite{SD}. Thus, the set
$\Delta^{\ir}_{\ge 0}$ is completely determined by the set
$\NEF(S)$. In what follows, we will use the set $\NEF(S)$ instead of
$\Delta^{\ir}_{\ge 0}$ since it is more convenient to work with.

Since $S$ is hyperbolic, the cone
$$
V(S\otimes \br)=\{x \in S\otimes \br\ |\ x^2>0 \}
$$
is the union of two half cones $\pm V^+(S\otimes \br)$ where
$V^+(S\otimes \br)$ contains the class of a hyperplane section.
It is easy to see that
$$
\NEF(S)=\{h \in S\ |\ h \in \overline{V^+(S\otimes \br)}-\{0\}\ \text{and}\
h\cdot \Delta_{-2}^{\ir}\ge 0\}.
$$
Thus, $\NEF(S)$ is completely defined by $\Delta_{-2}^{\ir}$.
Moreover, there exists a group-theoretical description of both sets.
Let
$$
\br_{++}\M=\{x \in \overline{V^+(S\otimes \br)}-\{0\}\ |
\ x\cdot \Delta_{-2}^{\ir}\ge 0\}
$$
be a cone and $\M=\br_{++}\M/\br_{++}$ its set of rays. Then
$\NEF(S)=S \cap \br_{++}\M$.
Let $W^{(2)}(S)\subset O(S)$ be the group generated by all reflections
$s_\delta: x \mapsto x+(x\cdot \delta)\delta$, $x \in S$ of the lattice
$S$ in elements $\delta \in S$ with $\delta^2=-2$. It is easy to see
that the group $W^{(2)}(S)$ is discrete in the corresponding
hyperbolic space $\La^+(S)=V^+(S)/\br_{++}$ and $\M$ is the
fundamental domain of $W^{(2)}(S)$ with the set
$\Delta^{ir}_{-2}$ of vectors orthogonal to $\M$. It means that
$\delta \in \Delta^{\ir}_{-2}$ if and only if $\delta \in S$,
$\delta^2=-2$ and the inequality
$\delta\cdot x\ge 0$ defines a face of $\M$ (or
of $\br_{++}\M$) of codimension one. This gives
the description of the both sets
$\NEF (S)$ and $\Delta^{\ir}_{-2}$ of $S$ in terms of the group
$W^{(2)}(S)$: the real convex
cone $\br_{++}\NEF(S)$ is a fundamental domain for the
group $W^{(2)}(S)$ acting in $V^+(S\otimes \br)$ with
the set of orthogonal vectors  $\Delta^{\ir}_{-2}$.

Let $\EF(S)$ be the set of all effective elements of $S$ and
$\EF(S)_{\ge -2}$, $\EF(S)_{-2}$ and $\EF(S)_{\ge 0}$
are the sets of all elements $x\in \EF(S)$ with $x^2 \ge -2$, $x^2=-2$ and
$x^2 \ge 0$ respectively. Using Riemann-Roch Theorem, one can see that
$$
\EF(S)_{-2}=\{\delta \in S\ |\ \delta^2=-2,\ \
\delta\cdot \NEF(S) \ge 0\},\ \
\EF(S)_{\ge 0}=S\cap \overline{V^+(S\otimes \br)}-\{0\},
$$
and \ $\EF(S)_{\ge -2}=\EF(S)_{-2}\cup \EF(S)_{\ge 0}$.

\subhead
2. Kac--Moody algebras associated to a K3 surface
\endsubhead

In this section, we define Kac--Moody algebras associated with a
K3 surface $X$ with the Picard lattice $S=S_X$.
(See \cite{Ka1}, \cite{Ka2}, \cite{Ka3}, \cite{Bo1} and \cite{GN}
for the theory of Kac--Moody algebras.) An algebra will be a
generalized Kac--Moody (Lie) superalgebra without odd real
simple roots. It is defined by a set ${}_s\Delta$ of
simple roots which is divided in a set of simple real (even)
roots ${}_s\Delta^{\re}$ and a set of simple imaginary roots
${}_s\Delta^{\im}$; the last set is divided in a set of even simple
imaginary roots ${}_s\Delta_{\overline 0}^{\im}$ and
a set of simple odd imaginary roots ${}_s\Delta_{\overline 1}^{\im}$.

We put
${}_s\Delta^{\re}=\Delta^{\ir}_{-2}$ and the sets
${}_s\Delta_{\overline 0}^{\im}$,
${}_s\Delta_{\overline 1}^{\im}$ are some sequences
of nef elements of $S$. Each imaginary root $\alpha$ defines
an element of $\NEF(S)$ but one can repeat each element of $\NEF(S)$
finite number of times in each set ${}_s\Delta_{\overline 0}^{\im}$ and
${}_s\Delta_{\overline 1}^{\im}$.

The {\it generalized Kac--Moody superalgebra without
odd real simple roots}
$\geg=\geg^{\prime\prime}(S, $ ${}_s\Delta^{\im})$
is a Lie superalgebra generated by
$h_r, e_r, f_r$ where $r \in {}_s\Delta$. All $h_r$ are even,
$e_r, f_r$ are even (respectively odd) if
$r$ is even (respectively odd).
The algebra has the following defining relations:

\vskip5pt

(1) The map $r \mapsto h_r$ for $r\in {}_s\Delta$ gives an embedding
of $S\otimes \br$ into $\geg^{\prime\prime}(S, {}_s\Delta)$ as
an abelian subalgebra (it is even since all $h_r$ are even).
In particular, all elements $h_r$ commute.

(2) $[h_r, e_{r^\prime}]=-(r\cdot r^\prime)e_{r^\prime}$, and \
$[h_r, f_{r^\prime}]=(r\cdot {r^\prime})f_{r^\prime}$.

(3) $[e_r, f_{r^\prime}]=h_r$ if $r=r^\prime$, and is \ $0$\  if
$r \not=r^\prime$.

(4) $(\text{ad\ } e_r)^{1+r\cdot r^\prime}e_{r^\prime }=
(\text{ad\ } f_r)^{1+r\cdot r^\prime}f_{r^\prime }=0$
if $r\not= r^\prime$ and $r \in {}_s\Delta^{\re}$.

(5) If $r\cdot r^\prime=0$, then $[e_r, e_{r^\prime}]=[f_r,f_{r^\prime}]=0$.

\vskip5pt

The superalgebra $\geg=\geg ^{\prime\prime}(S, {}_s\Delta^{\im})$
is graded by $S$ as follows. Let
$$
\widetilde{Q}_+ =\sum_{\alpha \in {}_s\Delta}{\bz_+\alpha}\subset S
$$
be the integral cone (semi-group) generated by all simple roots.
We have
$$
\geg=\left( \bigoplus_{\alpha \in \widetilde{Q}_+} {\geg_\alpha} \right)
\bigoplus \left(S\otimes \br \right) \bigoplus
\left( \bigoplus_{\alpha \in -\widetilde{Q}_+} {\geg_\alpha} \right)
$$
where $e_r$ and $f_r$ have degree $r\in \widetilde{Q}_+$ and
$-r \in -\widetilde{Q}_+$ respectively ($r \in {}_s\Delta$);
and $\geg_0=S\otimes \br$.
A non-zero $\alpha \in \pm \widetilde{Q}_+$ is called  a
{\it root} if $\geg_\alpha$ is
non-zero. Let $\Delta$ be the set of all roots and
$\Delta_{\pm}=\Delta\cap \pm \widetilde Q_+$.
For a root $\alpha \in \Delta$ we set
$\mult_\0o \alpha =\dim \geg_{\alpha,\0o}$,
$\mult_\1o \alpha = -\dim \geg_{\alpha,\1o}$ and
$$
\mult~\alpha=\mult_\0o \alpha +\mult_\1o \alpha=
\dim \geg_{\alpha ,\0o} - \dim \geg_{\alpha,\1o}.
$$
The $\mult~\alpha$ is called the multiplicity of $\alpha$.
According to the general theory of Kac--Moody algebras, the
set of roots is the union of real and imaginary roots:
$\Delta=\Delta^{\re}\cup \Delta^{\im}$.
The set of real roots is $\Delta^{\re}=W^{(2)}(S)({}_s\Delta^{\re})$
(in particular, $\alpha^2=-2$ if $\alpha \in \Delta^{\re}$).
The set of imaginary roots is
$\Delta^{\im}=\{\alpha \in \Delta\ |\ \alpha^2\ge 0\}$.
It follows that
$\Delta^{\re}_+:=\Delta^{\re}\cap \Delta_+=\EF_{-2}$ and
$\Delta^{\im}_+:=\Delta^{\im}\cap \Delta_+\subset \EF_{\ge 0}$.
If $\alpha \in \Delta^{\re}$, then $\mult_\0o \alpha =1$,
$\mult_\1o \alpha =0$ and
$\mult~\alpha=1$. Thus, we can rewrite
the decomposition above using ``geometry" of K3 as follows:
$$
\geg=\left( \bigoplus_{\alpha \in \EF(S)_{\ge -2}} {\geg_\alpha} \right)
\bigoplus \left(S\otimes \br \right) \bigoplus
\left( \bigoplus_{\alpha \in -\EF(S)_{\ge -2}} {\geg_\alpha} \right).
$$
Here $\geg_\alpha =0$ if $\alpha\not\in \Delta$.
\smallpagebreak

In what follows, we restrict ourselves
considering $S$ with a lattice Weyl vector.

\definition{Definition} An element $\rho \in S\otimes \bq$ is
called a {\it lattice
Weyl vector} if $\rho\cdot \delta =1$ for any
$\delta \in {}_s\Delta^{\re}=\Delta^{\ir}_{-2}$.
\enddefinition

There are three cases when a lattice Weyl vector does exist:

\smallpagebreak

(i) $\Delta^{\ir}_{-2}=\emptyset$,
then we can take any $\rho \in S\otimes \bq$;

\smallpagebreak

(ii) $\dim S=2$ and  $\Delta^{\ir}_{-2} \not=\emptyset$,
then the set $\Delta^{\ir}_{-2}$ is linearly independent and does not
contain more than $2$ elements;

\smallpagebreak

(iii) $\dim S\ge 3$ and  $\Delta^{\ir}_{-2}\not=\emptyset$ and the
lattice Weyl vector $\rho$ exists. It follows
from general results \cite{N4}, \cite{N5} and \cite{N10}, that the
set of hyperbolic lattices $S$ with this property is finite up to
isomorphism. These lattices $S$ are divided in two classes.
Firstly, it is easy to see that
$\rho$ is a nef element of $S$ and for large
$n \in \bn$ the linear system $|h|$ of $h=n\rho \in S$ is free.
If $\rho^2>0$ (this case is called {\it elliptic}),
the linear system $|h|$ gives an embedding of $X$
into a projective space such that all non-singular rational curves on $X$
have the same degree $n$. All these cases are known (see \cite{N3},
\cite{N7}, \cite{N8}). If $\rho^2=0$ (this case is called {\it parabolic}),
then $|h|$ gives an elliptic fibration of $X$ over a projective line
such that all non-singular rational curves of $X$ have the same degree $n$
over the projective line. The list of such $S$ is not known yet.

We would like to mention that in the case (iii) the fundamental polyhedron
$\M=\br_{++}\NEF(S)/\br_{++}$ for the action of
$W^{(2)}(S)$ in the hyperbolic space $\La ^+(S)$ is a very right
and beautiful polyhedron: it is a fundamental polyhedron for a reflection
group and it is touching of a sphere with the center $\br_{++}\rho$.
\smallpagebreak

\smallpagebreak

The case (iii) is especially interesting for us because it is very
exceptional: there is only finite number of possibilities. Moreover, we
want to get some relations between the sets $\NEF(S)$ and $\EF(S)$
which are different only for the cases (ii) and (iii). The case (iii)
is also related with multi-dimensional ($\dim \ge 3$) automorphic forms.

Later on we assume that $S$ has a lattice Weyl vector $\rho$. For
$a \in \NEF(S)$, let $m(a)_\0o^\prime $, $m(a)_\1o^\prime$
are equal to the numbers
of times we repeat $a$ in the sequences ${}_s\Delta_{\overline 0}^{\im}$
and   ${}_s\Delta_{\overline 1}^{\im}$ respectively. We set
$m(a)^\prime =m(a)_\0o^\prime - m(a)_\1o^\prime $. Let $a_0$ be
a primitive element of $\NEF(S)$ with $a_0^2=0$. In this case
we define ``corrected" invariants $m(ta_0)$, $t \in \bn$,
using the identity of power series:
$$
\prod_{n \in \bn}(1-q^n)^{m(na_0)^\prime }=
1-\sum_{t \in \bn}{m(ta_0) q^t} .
$$
For $a \in \NEF(S)$ with $a^2>0$ we set $m(a)=m(a)^\prime$.

We have the following Weyl--Kac--Borcherds denominator
identity for Kac--\-Moody superalgebra
$\geg=\geg^{\prime\prime}(S, {}_s\Delta)$ (see \cite{Ka1}, \cite{Bo1}
and \cite{GN}):
$$
\multline
\Phi(z):= \hskip-10pt \sum_{w\in W^{(2)}(S)}{\hskip-15pt \det(w)
\left(\exp(2\pi i (w(\rho)\cdot z))
- \hskip-10pt \sum_{ a \in \NEF(S)}\hskip-10pt
{m(a)\exp(2\pi i (w(\rho+a)\cdot z))}\right)}\\
\vspace{2\jot}
=\exp{\left(2\pi i(\rho\cdot z)\right)}
\prod_{\alpha \in \EF(S)_{\ge -2}}
{\left( 1-\exp{ \left(2\pi i (\alpha \cdot z)\right)}\right)^{\mult~\alpha}},
\endmultline
\tag*
$$
where $z$ belongs to the complexified cone
$\Omega(S)=S\otimes \br +iV^+(S\otimes \br)$ of $V^+(S\otimes \br)$.
The function $\Phi (z)$ is called the {\it denominator function}
of $\geg (S, {}_s\Delta^{\im})$.

Considering different sequences $\Delta^{\im}$ of
imaginary roots from $\NEF(S)$ we get different denominator
identities which one can consider as multi-dimensional identities
relating the sets of effective and nef (or irreducible) elements of $S$.
Actually these identities depend only on the integral function
$m(a)$, $a \in \NEF(S)$. If this function is given, one can calculate
$m(a)^\prime$ and find all possible non-negative integers
$m(a)_\0o^\prime$,  $m(a)_\1o^\prime$ with
$m(a)^\prime =m(a)_\0o^\prime - m(a)_\1o^\prime$
which define Kac--Moody superalgebras
$\geg=\geg^{\prime\prime}(S, {}_s\Delta)$ with the fixed denominator function.
One also can consider the function \thetag* as some kind of integral
$$
\Phi (z)= \int_{C \subset X}{\xi (C, z)}
$$
along effective curves on the K3 surface $X$ with $S_X=S$. This integral
could be correctly defined because we only use effective
classes in the formula
\thetag{*}.

\subhead
3. A variant of Mirror Conjecture
\endsubhead

Now we consider the hyperbolic lattice $S$ using the model (B).
We consider only
the simplest case when
$$
T=S\oplus U(k),\ \  k\in \bn,\ \
U(k)=\left(
\matrix
0&k\\
k&0
\endmatrix \right).
$$
Let $c_1,c_2$ be the bases of $U(k)$ with this intersection matrix.
Then
$z\mapsto \bc(z\oplus (-z^2/2)c_1\oplus (1/k)c_2)$ defines an
embedding corresponding to the cusp defined by $c_1$
of the complexified cone
$\Omega(S)$
to the connected component $\Omega(T)_0$ of the complex
domain of type IV
$$
\Omega(T)=\{\bc\omega \subset T\otimes \bc\ |\ \omega^2=0,\
\omega\cdot \overline{\omega}>0 \}.
$$
The choice of $\omega_0=z\oplus (-z^2/2)c_1\oplus (1/k)c_2\in \bc\omega_0$
is determined by the normalization
$\omega_0\cdot c_1=1$. For this normalization,
the local moduli of K3 are identified with $S\otimes \bc$ and Yukawa
coupling coincides with the intersection pairing of the lattice $S$.
This normalization is prescribed by mirror symmetry for K3.
The quotient $(\M_{T^\perp})_0=O(T)_+\setminus \Omega(T)_0$ is a connected
component of the dual (mirror) family of K3 surfaces
for an appropriate
subgroup $O(T)_+$ of finite index of $O(T)$.

\proclaim{Mirror Conjecture} There exists a choice of $k\in \bn$ and
a sequence
${}_s\Delta^{\im}\subset \NEF(S)$ of imaginary roots such that
the denominator function $\Phi (z)$ of
$\geg^{\prime\prime} (S, {}_s\Delta^{\im})$
is a holomorphic automorphic form with respect to $O(T)_+$ on the
domain $\Omega (T)_0$, i.e. $\Phi (z)$ is an
``algebraic function" on the dual moduli $\M_{T^\perp}$ (model (B)).
The form
 $\Phi(z)$ has the following sense from the point of view
of the model (A):\ \ $\Phi (z)$ is written in the form \thetag*
using  ``geometry of curves" (effective and irreducible or nef classes
of divisors) of a general member $X$ with $S_X=S$
of the family $\M_S$ and it gives an identity
\thetag{*} between effective and nef divisor classes on $X$.

Moreover, we suppose that for the automorphic form
$\Phi(z)$ it is possible
to give exact formulae for Fourier coefficients $m(a)$ of the
left side and multiplicities $\mult~\alpha$ of the right side of
\thetag*. Besides, the generalized Kac--Moody superalgebra
$\geg^{\prime\prime} (S, {}_s\Delta^{\im})$
should also be related with geometry of
curves and moduli of K3 (and conjecturally with some physics).
\endproclaim

\smallpagebreak

It is very important that the zero divisor of $\Phi (z)$ in the
domain where the product \thetag{*} converges has multiplicity one
and is contained in the discriminant
$$
{\Cal D}=O(T)_+\setminus (\bigcup_{\delta \in T, \delta^2=-2}{D_\delta})
$$
of moduli $\M_{T^\perp}$ where
$D_\delta=\{\bc \omega \in \Omega (T)\ | \ \omega\cdot \delta=0\}$.
Therefore,
in some sense, $\Phi (z)$ shows how far we are from the discriminant.

\smallpagebreak

In the rest part of the paper we give several examples
when this conjecture is valid.

\subhead
4. Example 1
\endsubhead

For the first example, $\dim S=3$, $S\cong 2U(4)\oplus
\langle -2 \rangle$. (In what follows we denote by $K(t)$ a lattice
which one gets by multiplication on $t \in \bq$ of the form of the
lattice $K$.)
The set $\Delta^{\ir}_{-2}(S)=\{\delta_1, \delta_2, \delta_3\}$
generates the lattice $S$
and has
the intersection matrix
$$
(\delta_i\cdot \delta_j)=
\left(
\matrix
 -2 & \hphantom{-}2 & \hphantom{-}2 \\
 \hphantom{-} 2 &-2 & \hphantom{-}2 \\
  \hphantom{-}2 & \hphantom{-}2 & -2
\endmatrix
\right)
$$
which defines the lattice $S$.
The fundamental polyhedron
$\M$ is the right triangle with the vertices at infinity. The lattice
Weyl vector $\rho$ is equal to $\rho=(\delta_1+\delta_2+\delta_3)/2$.
The element $h=2\rho$ has the square
$h^2=6$ and the linear system $|h|$ gives an embedding of a K3 surface
$X$ with $S_X=S$ as an intersection of a quadric and a cubic in $\Bbb P^4$.
For this embedding, all non-singular rational curves on $X$ are three conics
corresponding to $\delta_1, \delta_2, \delta_3$. Their sum is a hyperplane
section of $X$. The lattice $T$ is equal to
$T=U(4) \oplus S\cong 2U(4)\oplus \langle -2 \rangle$.
The orthogonal complement $T^\perp$ is isomorphic to a hyperbolic lattice
$S^\prime\cong U(4)\oplus K$ where $K$ is a negative definite lattice
of rank $15$ with
the discriminant quadratic form $q_{U(4)}\oplus q_{\langle 2 \rangle}$.
It follows from results of \cite{N2} that the lattice $S^\prime$ is
unique and the moduli space of K3 surfaces
$\M_{S^\prime}$ is irreducible. (It would be
very interesting to determine this family using
equations and to give an algebraic description of the automorphic
form $F_1(Z)$ which we shall describe. We hope to do this later.)

\smallpagebreak

Let us consider another bases
$f_2, f_3, f_{-2}$ of $S\otimes \bq$ where
$$
\delta_1 =2f_2-f_3,\  \delta_2=2f_{-2}-f_3,\ \
\delta_3=f_3.
$$
These elements have the intersection matrix
$$
(f_i\cdot f_j)=
\left(\matrix
0 & \hphantom{-}0 & 1\\
0 & -2 & 0\\
1  & \hphantom{-}0 & 0
\endmatrix
\right).
$$
Thus the lattice $S$ is a finite index sublattice of
$M_0=\bz f_2\oplus \bz f_3 \oplus \bz f_{-2}$. We have
$M_0=U\oplus \langle -2 \rangle$ where $U=\bz f_2 \oplus \bz f_{-2}$ and
$\bz f_3=\langle -2 \rangle$. These lattices are
related as follows:
$S=2(M_0)^\ast$.
We consider coordinates $(z_3,z_2,z_1)$
where $z=z_3f_2+z_2f_3+z_1f_{-2}\in M_0\otimes \bc=S\otimes \bc$.
We introduce the lattice $L=U\oplus M_0$ where
$U=\bz f_1\oplus \bz f_{-1}$ with $f_1^2=f_{-1}^2=0$ and $f_1\cdot f_{-1}=1$.
We use $z$ as a coordinate for the point
$Z = \bc ((-z^2/2)f_1+z+f_{-1})\in \Omega(L)$ of the domain
$\Omega (L)$ of the type IV corresponding to $L$.
We also identify $z$ with the matrix
$$
Z=\left(
\matrix
z_1&z_2\\
z_2&z_3
\endmatrix
\right) \in \bh_2
$$
where $\bh_2$ is the Siegel upper-half plane.

\smallpagebreak

Let us consider the classical function $\Delta_5(Z)$ (see \cite{F})
which is the product of all even
theta-constants
$$
\Delta_5 (Z)=\prod_{(a,b)}\vartheta_{a,b}(Z), \qquad
(Z=\pmatrix z_1&z_2\\z_2&z_3\endpmatrix\in \bh_2)
$$
with
$$
\vartheta_{a,b}(Z)=\sum_{l\in \bz^2}
\exp{\bigr(\pi i (Z[l+\frac 1{2}a]+{}^tbl)\bigl)}
\qquad\qquad (Z[l]={}^tlZl).
$$
The product is taken over all vectors $a,b\in (\bz/2\bz)^2$
such that ${}^t ab\equiv 0\mod 2$.
(There are exactly ten different $(a,b)$.)
This is the automorphic form of weight $5$ with a character
with respect to $Sp_4(\bz)/\{\pm E_4\}\cong O^+(L)/\{\pm E_5\}$
where $O^+(L)$ is the subgroup of $O(L)$ which fixes two connected
components of $\Omega(L)$ (see \cite{GN}).
The function
$F_1(Z)=\frac{1}{64}\Delta_5 (Z)$ has integral coefficients .

\proclaim{Theorem 1} The function $F_1(Z)$ gives the solution of
Mirror Conjecture of Sect.3 for the lattice $S$ and
$U(4)=\bz c_1\oplus \bz c_2$ where
$c_1=2f_1, c_2=2f_{-1}$. Therefore
$F_1(Z)$ is an ``algebraic function" on the moduli $\M_{T^\perp}$
where $T=U(4)\oplus S$, and it defines
an identity \thetag{*} for $S=S_X$ of the general member $X$ of
the mirror family $\M_S$. Moreover, it defines the corresponding Kac--Moody
superalgebras $\geg^{\prime\prime}(S, {}_s\Delta^{\im})$.
\endproclaim

\demo{Proof}  The function $F_1(Z)$ as a function
on $\Omega (L)$ is automorphic with respect to $O^+(L)$. We have the
equality
$T=2L^\ast$ because $U(4)=2U^\ast$ and $S=2(M_0)^\ast$. It follows
that $F_1(Z)$ is automorphic with respect to $O^+(T)=O^+(L)$ and
defines then an ``algebraic function" on the moduli
$\M_{T^\perp}=O^+(T)\setminus \Omega (T)$.

It is proved in \cite{GN} that for the coordinate $z$ which we
introduced above, the function $F_1(Z)$ can be written in the
form
$$
\split
F_1(Z)= \hskip-10pt \sum_{w\in W^{(2)}(S)}{\hskip-5pt \det(w)
\left(\exp(\pi i (w(\rho)\cdot z))
- \hskip-10pt \sum_{ a \in \NEF(S)}\hskip-10pt
{m(a) \exp(\pi i (w(\rho+a)\cdot z))}\right)}\\
\vspace{2\jot}
=\exp{\left(\pi i(\rho \cdot z)\right)}
\prod_{\alpha \in \EF(S)_{\ge -2}}
{\left( 1-\exp{ \left(\pi i (\alpha \cdot z)\right)}\right)^{\mult~\alpha}}
\endsplit
\tag4.1
$$
with integral coefficients $m(a)$ and $\mult~\alpha$.
For $c_1=2f_1$ and $k=4$, we should consider the coordinate
$z^\prime =z/2$ (mirror symmetry coordinate) instead of the coordinate
$z$. For this coordinate $z^\prime$, from \thetag{4.1} we get
$$
\split
F_1(Z)= \hskip-10pt \sum_{w\in W^{(2)}(S)}{\hskip-5pt \det(w)
\left(\exp(2\pi i (w(\rho)\cdot z^\prime ))
- \hskip-10pt \sum_{ a \in \NEF(S)}\hskip-10pt
{m(a) \exp(2\pi i (w(\rho+a)\cdot z^\prime ))}\right)}\\
\vspace{2\jot}
=\exp{\left(2\pi i(\rho\cdot z^\prime )\right)}
\prod_{\alpha \in \EF(S)_{\ge -2}}
{\left( 1-\exp{\left(2\pi i (\alpha \cdot z^\prime)\right)}
\right)^{\mult~\alpha}}.
\endsplit
\tag4.2
$$
This proves Theorem 1.
\enddemo

To calculate the coefficients $m(a)$ and $\mult~\alpha$ in the
product formula \thetag{4.1}, we need two types of Jacobi
forms. The Jacobi form of the first type is the form of weight
$5$ and index $1/2$
$$
\psi_{5,\frac 1{2}}(z_1,z_2)
=\eta(z_1)^9\,\vartheta_{11}(z_1,\,z_2).
$$
Here
$\eta(z_1)=q^{\frac 1{24}}\,\prod_{n\ge 1}\ (1-q^n)$
is Dedekind eta-function
and
$$
\align
\vartheta_{11}(z_1,\,z_2)
&=\sum\Sb n\in \bz\endSb
\,(-1)^{n}
\exp{\bigl(\frac{\pi i\, }{4}(2n+1)^2 z_1+\pi i\,(2n+1)z_2\bigr)}\\
{}&=-q^{1/8}r^{-1/2}
\prod_{n\ge 1}\,(1-q^{n-1} r)(1-q^n r^{-1})(1-q^n)
\endalign
$$
is the  classical Jacobi theta-function,
where we put
$$
\split
z_1\in \bh_1=\{z_1=x+iy\in \bc\,|\,y>0\},\  z_2\in \bc,\\
\ q=\exp{ (2\pi i z_1)}, \ r=\exp{ (2\pi i z_2)},\ p=\exp{(2\pi i z_3)}.
\endsplit
$$
The holomorphic function
$\psi_{5,\frac 1{2}}(z_1,z_2)$ is a Jacobi form of index one-half
with a multiplier system. It means that the following identities
are satisfied
$$
\align
\psi_{5,\frac 1{2}}(\frac{az_1+b}{cz_1+d},\,\frac{z_2}{cz_1+d})&=
v^{12}_\eta(g)(cz_1+d)^{5}
\exp{(\pi i\,\frac{ cz_2^2}{cz_1+d})}\,
\psi_{5,\frac 1{2}}(z_1 ,\,z_2 ),
\\
\vspace{2\jot}
\psi_{5,\frac 1{2}}(z_1, z_2+pz_1+q)&=
(-1)^{p+q}\exp{(-\pi i (p^2z_1 +2p z_2))}\,\psi_{5,\frac 1{2}}(z_1, z_2),
\endalign
$$
where $p,q\in \bz$ and
$g=\left(\smallmatrix a&b\\c&d\endsmallmatrix\right) \in SL_2(\bz)$
and
$$
\eta(\frac{a\tau+b}{c\tau+d})=v_\eta(g)(c\tau+d)^{1/2}\eta(\tau).
$$
Here $v_\eta(g)$ is a $24$th root of unity.

Let us consider the Fourier coefficients of $\eta^d(\tau)$
$$
q^{\frac{d}{24}}\prod_{n\ge 1}\ (1-q^n)^d=\sum_{m} {\tau_d(m)q^m}.
$$
Then we have  the  following Fourier expansion
of $\psi_{5,\frac{1}{2}}(z_1,z_2)$:
$$\align
\psi_{5,\frac{1}{2}}(z_1,z_2)&=\eta(z_1)^9\,\vartheta_{11}(z_1,\,z_2)\\
&=\sum\Sb n,l \equiv 1\, mod\,2 \\ n>0,\,4n-l^2>0\endSb
(-1)^{\frac{l-1}2}\,\tau_9(4n-l^2)\,\exp{(\pi i(nz_1+lz_2))}.
\endalign
$$

The second type of Jacobi forms which we need
are special Jacobi forms of weight zero
(weak forms in terms of \cite{EZ}). The
ring of all weak Jacobi forms has two generators
as an algebra over $SL_2(\bz)$-modular forms (see \cite{EZ, \S 9}).
One of these generators is  the function
$$
\phi_{0,1}(z_1, z_2)=\frac 1{144 \Delta(z_1)}
\bigr(E_4^2(z_1)E_{4,1}(z_1,z_2)-E_6(z_1)E_{6,1}(z_1,z_2)\bigl)
$$
where
$$
\Delta(z_1)=q\prod_{n}(1-q^n)^{24},\
E_4(z_1)=1+240 \sum_{n\ge 1}\sigma_3(n)q^n, \
E_6(z_1)=1-504 \sum_{n\ge 1}\sigma_5(n)q^n
$$
are the cusp form of weight $12$,
the  Eisenstein series for $SL_2(\bz)$
and $E_{k,1}(z_1,z_2)$ is  the Jacobi-Eisenstein series
of  weight $k$  and index one. This series has the following
integral Fourier coefficients (see \cite{EZ, \S 2})
$$
E_{k,1}(z_1,z_2)=\zeta(3-2k)^{-1}\sum\Sb n,l\in \bz\\ 4n-l^2\ge 0\endSb
H(k-1,4n-l^2)  q^n r^l,
$$
where $H(k,N)=L_{N}(1-k)$ are H. Cohen's numbers (see \cite{C}).
We recall that $\zeta(-5)=-1/252$ and  $\zeta(-9)=1/132$.

The form $\phi_{0,1}$ has  the Fourier expansion with
integral coefficients
$$\align
\phi_{0,1}(z_1, z_2)&=
\sum\Sb n,l\in \bz, n\ge 0\\ 4n-l\ge -1\endSb
f(n,l)\,\exp{(2\pi i\,(nz_1+lz_2))}\\
\vspace{ 2\jot}
{}&=(r^{-1}+10+r)+q(10r^{-2}-
64r^{-1}+108-64r+10r^2)+\dots\
\endalign
$$
which depend only on the ``norm'' $4n-l^2$ of $(n,l)$
$$
f(n,l)=c_1(4n-l^2)\qquad\text{and}\qquad c_1(m)=0\ \text{ for }m<-1.
$$
Moreover, the function
$$
C_1(z_1)=\sum\Sb m\ge -1\\ m\equiv 0,3\, mod\, 4\endSb
c_1(m)q^m
$$
is an automorphic form of weight $-\frac{1}2$.
It is easy to get a formula for this function
using the Cohen's modular forms of half integral weight $k-\frac1{2}$
$$
{\Cal H}_{k-1}(z_1)=\sum_{n\ge 0} H(k-1,n)q^n
$$
(see \cite{C} and \cite{EZ, \S 5}).
One has
$$\align
C_1(z_1)&=\frac{1}{12 \Delta(4z_1)}
\bigl( 11 E_6(4z_1)\Cal H_5(z_1) -21 E_8(4z_1) \Cal H_3(z_1) \bigr)\\
{}&=q^3+10q^4+ \cdots .
\endalign
$$
We would like to note that  the function
$(11 E_6(4z_1)\Cal H_5(z_1) -21 E_8(4z_1) \Cal H_3(z_1))$
is the cusp form of weight $11\frac1{2}$ for $\Gamma_0(4)$.

Using the functions introduced above and Theorem 4.1
of  \cite{GN}, we can write the identity \thetag{4.1}
in the following form:
$$
\multline
F_1(Z)=\sum\Sb n,l,m \equiv 1\,mod\,2\\
\vspace{0.5\jot} n,m>0\endSb
-\sum_{d|(n,l,m)} (-1)^{\frac{l+d}2}\,
\tau_9\,(\frac{4nm-l^2}{d^2})\,
q^{n/2}\,r^{l/2}\,p^{m/2}
\\
\vspace{2\jot}
=(qrp)^{1/2} \prod
\Sb n,\,l,\,m\in \Bbb Z\\
\vspace{0.5\jot}
(n,l,m)>0\,\endSb
\bigl(1-q^n r^l p^m\bigr)^{c_1(4nm-l^2)},
\endmultline
$$
where $(n,l,m)>0$ means that $n\ge 0$, $m\ge 0$,
$l$ is an  arbitrary integral if $n>0$ or $m>0$ and
$l<0$ if $n=m=0$.

\subhead
5. Example 2
\endsubhead

For this example, $\dim S=3$, $S\cong 2U(8)\oplus \langle -2 \rangle$.
The set $\Delta^{\ir}_{-2}(S)=\{e_1, e_2, e_3, e_4\}$
generates the lattice $S$
and has
the intersection matrix
$$
e_i\cdot e_j=
\left(
\matrix
-2 &  \hphantom{-}2 & \hphantom{-}6 & \hphantom{-}2 \\
 \hphantom{-}2 & -2 & \hphantom{-}2 & \hphantom{-}6\\
\hphantom{-}6 &  \hphantom{-}2 &-2 & \hphantom{-}2 \\
 \hphantom{-}2 &  \hphantom{-}6 & \hphantom{-}2 &-2
\endmatrix
\right)
$$
which defines the lattice $S$.
The fundamental polyhedron
$\M$ is the right quadrangle with the vertices at infinity. The lattice
Weyl vector $\rho$ is given by the equality
$\rho =(e_1+e_3)/4=(e_2+e_4)/4$.
The element $h=4\rho$ has the square
$h^2=8$ and the linear system $|h|$ gives an embedding of a K3 surface
$X$ with $S_X=S$ as an intersection of three quadrics in $\Bbb P^5$ (this
follows easily from general results of \cite{SD}).
For this embedding, all four
non-singular rational curves on $X$ have degree $4$.
The curves $e_1+e_3$ and $e_2+e_4$ give
two hyperplane sections of $X$.
The lattice $T$ is
$T=U(8)\oplus S\cong 2U(8)\oplus \langle -2 \rangle$.
The orthogonal complement $T^\perp$ is isomorphic to a hyperbolic lattice
$S^\prime\cong U(8)\oplus K$ where $K$ is a negative definite lattice of
the rank 15 with
the discriminant quadratic form $q_{U(8)}\oplus q_{\langle 2 \rangle}$.
It follows from
results of \cite{N2}, that the lattice $S^\prime$ is unique and
the moduli space $\M_{S^\prime}$ is irreducible
\smallpagebreak

We describe below an automorphic form for this case
which gives the solution of the mirror conjecture in Sect.3 for this case.
We consider a hyperbolic lattice $M_0$ with the
bases $f_2, f_3, f_{-2}$ and the intersection matrix
$$
(f_i\cdot f_j)=
\left(\matrix
0 & \hphantom{-}0 & 1\\
0 & -4 & 0\\
1  & \hphantom{-}0 & 0
\endmatrix
\right).
$$
Let us take the hyperbolic plane $U$ with the standard bases
$f_1$, $f_{-1}$ where $f_1^2=f_{-1}^2=0$, $f_1\cdot f_{-1}=1$, the
lattice $L=U\oplus M_0$ and the domain $\Omega (L)$. We use
the coordinate
$z^\prime=z_3^\prime f_2+z_2^\prime f_3+z_1^\prime f_{-2}$ for a point
$Z^\prime = \bc((-(z^\prime)^2/2)f_1+z^\prime + f_{-1})$ of this domain.
We define
$$
\{\delta_1=-f_3,\ \delta_2=4f_2+f_3,\
\delta_3=4f_2+3f_3+4f_{-2},\ \delta_4=f_3+4f_{-2}\}\subset M_0.
$$
Then $\delta_i\cdot \delta_j=2e_i\cdot e_j$. Thus,
the sublattice $M_{II}\subset M_0$ generated by
$\delta_1,...,\delta_4$ is isomorphic to $S(2)\cong U(16)\oplus
\langle -4 \rangle $. Equivalently,
$S=M_{II}(1/2)$. We identify this lattices replacing $e_i$ by $\delta_i$.

In \cite{GN, \S 5} there was constructed
an automorphic cusp form $F_2(Z^\prime)$ of weight $2$ with a character
with respect to $O^+(L)/\{\pm E_5\}$. This function has the
following representation with integral coefficients
$$
F_2(Z^\prime)=
\sum_{w\in W^{(2)}(S)}{ \hskip-10pt \det(w)
\left(\exp(\frac{\pi i}{2} (w(\rho)\cdot z^\prime))
- \hskip-10pt  \sum_{ a \in \NEF(S)}
{\hskip-10pt m(a)\exp(\frac{\pi i}{2} (w(\rho+a)\cdot z^\prime ))}\right)}
$$
$$
=\exp{\left(\frac{\pi i}{2}(\rho\cdot z^\prime )\right)}
\prod_{\alpha \in \EF(S)_{\ge -2}}
{\left( 1-\exp{ \left(\frac{\pi i}{2}
(\alpha\cdot z^\prime )\right)}\right)^{\mult~\alpha}}.
\tag5.1
$$

\proclaim{Theorem 2} The function $F_2(Z^\prime )$ gives the solution of
Mirror Conjecture of Sect.3 for the lattice
$S=M_{II}(1/2)\cong 2U(8)\oplus \langle -2 \rangle $ and
$U(8)=[\bz c_1\oplus \bz c_2](1/2)$ where
$c_1=4f_1, c_2=4f_{-1}$. Therefore,
$F_2(Z)$ is an ``algebraic function" on the moduli $\M_{T^\perp}$ (model
(B)) where $T=U(8)\oplus S$, and it defines
an identity \thetag{*} for $S=S_X$ of the general member $X$ of
the mirror family $\M_S$ (model (A)).
Moreover, it defines the corresponding Kac--Moody
superalgebras $\geg ^{\prime\prime} (S, {}_s\Delta^{\im})$.
\endproclaim

\demo{Proof}  The function $F_2(Z^\prime )$ as a function
on $\Omega (L)$ is automorphic with respect to $O^+(L)$. We have
$T(2)=4L^\ast$ because $U(16)=4U^\ast$ and $M_{II}=4(M_0)^\ast$.
It follows that $F_2(Z^\prime )$ is automorphic with respect to
$O^+(T)=O^+(T(2))=O^+(L)$ and
defines an ``algebraic function" on the moduli
$\M_{T^\perp}=O^+(T)\setminus \Omega (T)$.

For $c_1=4f_1, c_2=4f_{-1}$ and $U(8)$ we should use the
mirror symmetry coordinate
$z^{\prime\prime}=z^\prime /2$. Also
we should remember that $S=M_{II}(1/2)$.
{}From \thetag{5.1}, we get the identity
$$
F_2(Z^\prime)=\hskip-10pt
\sum_{w\in W^{(2)}(S)}{ \hskip-10pt \det(w)
\left(\exp(2\pi i(w(\rho)\cdot z^{\prime\prime}))
- \hskip-10pt  \sum_{ a \in NEF(S)}
{\hskip-10pt m(a)\exp(2\pi i (w(\rho+a)\cdot z^{\prime\prime} ))}\right)}
$$
$$
=\exp{\left(2 \pi i (\rho\cdot z^{\prime\prime})\right)}
\prod_{\alpha \in \EF(S)_{\ge -2}}
{\left( 1-\exp{ \left(2 \pi i
(\alpha\cdot z^{\prime\prime})\right)}\right)^{\mult~\alpha}}.
\tag{5.1$^\prime$}
$$
where we use the intersection pairing of the lattice $S$ and
$z^{\prime\prime} \in S\otimes \bc$.
This proves Theorem 2.
\enddemo

The function $F_2 (Z^\prime)$
(see \cite{GN, \S 5}) is connected
with the Jacobi functions
$$
\align
\phi_{0,2}(z_1,z_2)&=\frac 1{288\Delta_{12}(z_1)}
\bigl(E_4 (z_1)E_{4,1}^2(z_1,z_2)-E_{6,1}^2(z_1,z_2)\bigr)\\
\vspace{2\jot}
{}&=\sum\Sb n,l \endSb\  c_2 (8n-l^2)\,\exp{(2\pi i\, (n z_1 +l z_2))}
\endalign
$$
and
$$\align
\psi_{2,\frac{1}2}(z_1,z_2)&=-\eta^3(\tau)\vartheta_{11}(z_1,z_2)\\
{}&=\sum
\Sb n\equiv 1\,mod 4\,\ l \equiv 1\, mod\,2\endSb
(-1)^{\frac{l+1}2}\,\tau_3\,(2n-l^2)\,\exp{(\pi i(nz_1+lz_2))}.
\endalign
$$
The coefficients $\tau_3(n)$ are given by the Jacobi formula
$$
\eta^3(z_1)=\sum_{m\ge 1} \biggl(\frac {-4}{m}\biggr) mq^{m^2/4},
$$
where
$$
\biggl(\frac {-4}{m}\biggr)
=\cases \hphantom{-}1&\text{if }\ m\equiv \hphantom{-}1\mod 4\\
-1&\text{if }\ m\equiv -1\mod 4.
\endcases
$$
The numbers $c_2(n)$, which define the Fourier coefficients
of the Jacobi form $\psi_{2,\frac{1}{2}}$, are Fourier coefficients
of an automorphic  form  of weight $-1/2$.
One can express  them in terms of Cohen's numbers $H(3,N)$
and $H(5,N)$.
Using these functions we can rewrite the identity \thetag{5.1} in
the following form (see \cite{GN, \S 5})
$$\multline
F_2(Z^\prime )=\sum_{N^{}\ge 1}\ \sum\Sb 2mn-l^2=N^2\\
\vspace{0.5\jot}
 n,m\equiv 1\, mod\, 4\\
\vspace{0.5\jot}
n>0,\
l\equiv 1\, mod\,  2,\endSb
(-1)^{\frac{l+1}2}
\biggl(\frac {-4}{N}\biggr) N
\sum_{d\,|\,(n,l,m)} \biggl(\frac {-4}{d}\biggr)
\, q^{n/4} r^{l/2} p^{m/4}=\\
q^{1/4} r^{-1/2} p^{1/4}
\prod\Sb n,\,l,\,m \in \Bbb Z\\ (n,l,m)>0\endSb
(1- q^{n} r^{l} p^{m})^{c_2(8nm-l^2)},
\endmultline
$$
where $(n,l,m)>0$ means that $n\ge 0$, $m\ge 0$,
$l$ is an  arbitrary integral if $n+m>0$, and $l>0$ if $n=m=0$; \ \
 $p=\exp{(2\pi i z_1^\prime )}$,
$q=\exp{(2\pi i z_2^\prime )}$, $r=\exp{(2\pi i z_3^\prime )}$.

Let us consider the coordinate
$
Z=\left(
\matrix
z_1&z_2\\z_2&z_3
\endmatrix
\right)\in \bh_2
$
where $z_1=z_1^\prime$, $z_2=z_2^\prime$, $z_3=z_3^\prime /2$.
The function $F_2(Z)$ is the cusp form of weight 2 with
a character (with values in the group of fourth roots of unity)
for the double  extension of the paramodular group
$$
\Gamma_2:=
\left\{\pmatrix *&2*&*&*\\
                    *&*&*&2^{-1}*\\
                    *&2*&*&*\\
                    2*&2*&2*&* \endpmatrix
\in Sp_4(\Bbb Q),\quad\text{all } *\in \bz\right\}.
$$
This function is a lifting of Jacobi form $\psi_{2,\frac{1}{2}}$
(see \cite{G2}, \cite{G3}, \cite{G4}).

\subhead
6. Example 3
\endsubhead

For this example, $\dim S=10$, $S\cong U\oplus E_8(2)$. Let
$U=\bz c \oplus \bz e$ where $c^2=0$, $e^2=-2$ and $c\cdot e=1$.
Then
$$
\Delta^{\ir}_{-2}=\{ \delta \in S\ |\ \delta^2=-2,\ \delta\cdot c=1\}.
$$
For example, $e \in \Delta_{-2}^{\ir}$. This case is parabolic and
$\rho=c$. For a K3 surface $X$ with $S_X=S$, we have
$|\rho|:X \to {\Bbb P}^1$ is elliptic fibration. All non-singular
rational curves on $X$ are sections of this fibration.
Probably, this family of K3 surfaces had
first appeared in \cite{N3} (see also \cite{N8})
where $X \in \M_S$ were described as follows.
There exists an involution $\sigma$ on $X$ such that
$H^2(X, \bz)^\sigma=S$.
This involution is unique on $X$ and $\sigma^\ast\omega_X=-\omega_X$.
The set of points of $X$ fixed by this involution is union of two
non-singular fibers  (two elliptic curves)
of the fibration $|\rho|$ above. Let $Y$ be a
K3 surface with involution $\sigma$ on $Y$ such that the
set of points of $Y$ fixed by this
involution is union of two elliptic curves. Then
$S\cong H^2(Y, \bz)^\sigma \subset S_Y$ and $Y$
belongs to $\M_S$.

We can interpret results of \cite{Bo3} as construction of
a function $F_3(Z)$ which gives the
solution of Mirror Conjecture in Sect.3 for $U(2)$ ($k=2$).
Thus, for this case,
$T=U(2)\oplus S \cong U(2)\oplus U(1)\oplus E_8(2)$.
Then $S^\prime =T^\perp\cong U(2)\oplus E_8(2)$ and $\M_{S^\prime}$ is the
family of K3 surfaces which are universal coverings of Enriques surfaces
(``Enriques family"). In other words, $X \in M_{S^\prime}$ has an involution
$\sigma$ without fixed points. Then $X/\{1,\sigma\}$ is an Enriques surface.

\subhead
7. Questions
\endsubhead

It would be interesting to formulate Mirror Conjecture
of Sect.3 for hyperbolic lattices $S$ which do not have a lattice Weyl
vector for $\Delta^{\ir}_{-2}$. It is
certainly possible for some cases. Is this possible for arbitrary
hyperbolic lattices $S$?
What is an analog of Mirror Conjecture in Sect. 3 for Calabi--Yau 3-folds?

\enddocument

\end